\documentclass[a4paper]{ifacconf}

\usepackage[round]{natbib} 
\usepackage{graphicx}      
\usepackage{amsmath}
\usepackage{amssymb}
\usepackage{color}

\usepackage{float}
\graphicspath{ {pics/} } 
\usepackage{xparse}

\NewDocumentCommand{\INTERVALINNARDS}{ m m }{
	#1 {,} #2
}
\NewDocumentCommand{\interval}{ s m >{\SplitArgument{1}{,}}m m o }{
	\IfBooleanTF{#1}{
		\left#2 \INTERVALINNARDS #3 \right#4
	}{
		\IfValueTF{#5}{
			#5{#2} \INTERVALINNARDS #3 #5{#4}
		}{
			#2 \INTERVALINNARDS #3 #4
		}
	}
}
\begin{document}
\begin{frontmatter}

\title{\mbox{Boundary Observer for Congested Freeway }\protect\\   \mbox{\!\!\!\!\!\!\!\!\!\!\!\!\!\! Traffic State Estimation via Aw-Rascle-Zhang model}\protect}



\author[First]{Huan Yu} 
\author[Second]{Alexandre M. Bayen} 
\author[Third]{Miroslav Krstic}

\address[First]{University of California, San Diego, CA 92093 USA, \\(e-mail: huy015@ucsd.edu).}
\address[Second]{University of California, Berkeley, CA 94720-2284 USA, \\(e-mail: bayen@berkeley.edu)}
\address[Third]{University of California, San Diego, CA 92093 USA,\\ (e-mail: krstic@ucsd.edu)}

\begin{abstract}                
	This paper develops boundary observer for estimation of congested freeway traffic states based on Aw-Rascle-Zhang(ARZ) partial differential equations (PDE) model. Traffic state estimation refers to acquisition of traffic state information from partially observed traffic data. This problem is relevant for freeway due to its limited accessibility to real-time traffic information. We propose a boundary observer design so that estimates of aggregated traffic states in a freeway segment are obtained simply from boundary measurement of flow and velocity. The macroscopic traffic dynamics is represented by the ARZ model, consisting of $2 \times 2$ coupled nonlinear hyperbolic PDEs for traffic density and velocity. Analysis of the linearized ARZ model leads to the study of a hetero-directional hyperbolic PDE model for congested traffic regime. Using spatial transformation and PDE backstepping method, we construct a boundary observer with a copy of the nonlinear plant and output injection of boundary measurement errors. The output injection gains are designed for the error system of the linearized ARZ model so that the exponential stability of error system in the $L^2$ norm and finite-time convergence to zero are guaranteed. Simulations are conducted to validate the boundary observer design for nonlinear ARZ model without knowledge of initial conditions. 
\end{abstract}

\begin{keyword}
Aw-Rascle-Zhang model, boundary observer, traffic estimation, PDE backstepping method.
\end{keyword}

\end{frontmatter}

\section{Introduction}
   Traffic state estimation refers to foresee of traffic state information with a model by accessing partially observed traffic data and some prior knowledge of the traffic. 
   Traffic state estimation plays an important role in traffic management. In order to mitigate freeway traffic congestion, various control algorithms are developed for ramp metering or variable speed limit. The effective implementation of control algorithms on traffic infrastructure relies on accurate information of traffic state. Due to the financial and technical limitations, traffic state on freeways is difficult to measure everywhere. This topic has been widely studied and gained increasing attention in recent decades. 
   
   Freeway traffic dynamics in spatial and temporal domain are usually described with macroscopic models of the aggregated values of traffic states, including traffic density, velocity and flux.
   The Lighthill-Whitham-Richards (LWR) model proposed by \cite{LW} and \cite{R}, a first-order scalar hyperbolic PDE of density, is the most widely applied traffic models. Several studies have used LWR model for traffic states estimation in \cite{Cl}~\cite{Coi02}~\cite{Coi03}~\cite{kest} due to its simplicity. 
   
   This model fails to account for stop-and-go traffic, which does not obey the density-velocity relation in equilibrium. Second-order models consisting of nonlinear hyperbolic PDEs of traffic density and velocity, have been proposed to overcome the limitation of LWR model. The deviations from the equilibrium traffic relation are allowed in the second-order model. To estimate the non-equilibrium traffic states for congested traffic, the second-order models therefore need to be considered. 
   
The second-order Payne-Whitham (PW) model by \cite{p}~\cite{w} is used to develop an extended Kalman filter for state estimation in~\cite{Wa}. Compared with PW model, Aw-Rascle Zhang model by \cite{AW} and \cite{Zhang} improved the second-order model by successfully addressing anisotropic behavior of traffic and correcting the model's prediction of traffic waves. For this reason, ARZ model has been studied intensively for stop-and-go traffic by~\cite{Seibold1}~\cite{Seibold2}~\cite{Seibold3}~\cite{Bayen}~\cite{kerner}. A comprehensive review of different models and approaches in traffic estimation problem can be found in \cite{bayen}. Boundary observer for state estimation has never been addressed with ARZ PDE model.

Boundary control of the ARZ model has been studied through many recent efforts including~\cite{Bayen}~\cite{zpri}~\cite{KBM}~\cite{Huan1}~\cite{Huan2}. Boundary control and observer design using PDE backstepping method have been developed for $2\times 2$ coupled hyperbolic systems~\cite{rafael}\cite{Huan3} and the theoretical result for more  general hetero-directional hyperbolic systems developed in~\cite{long}~\cite{florent1}~\cite{florent3}.
The previous work by~\cite{Huan1}~\cite{Huan2} adopted the methodology for ARZ model. Boundary observers are designed in an effort to develop an output feedback for the linearized ARZ model. To the author's best knowledge, state estimation problem of nonlinear ARZ model has never been studied before.

This is the first result on boundary observer of nonlinear ARZ PDE model. Our main contribution in this work is that we generalize the previous results of linearized ARZ model and address the freeway traffic estimation problem from a more practical perspective. The boundary conditions are obtained from measurement of the field data. The current result can be easily represented and validated by real traffic data and thus paves the way for implementing this observer design in practice.

The outline of this paper is as follows: 	
we firstly introduce the nonlinear ARZ model, and analyze the linearized ARZ model for free and congested traffic. The boundary observer for linearized ARZ model is designed using the backstepping method and nonlinear boundary observer are developed using the output injections obtained from the linearized model. To validate our result, simulations of nonlinear ARZ PDE model and state estimation by nonlinear boundary observer are conducted. The estimation errors are analyzed.

\section{Problem statement}
We consider the traffic estimation problem for a stretch of freeway with length to be $L$. The macroscopic traffic dynamics is modeled by the ARZ model. We study the linearized ARZ model and discuss the characteristic speeds under free and congested traffic regime. 
\subsection{Aw-Rascle-Zhang Model}
The ARZ model for $(x,t)\in [0,L]\times \interval[{0, +\infty})$ is 

\begin{align}
\partial_t \rho + \partial_x( \rho v)=&0,\label{rho} \\
\partial_t v+(v-\rho p'(\rho))\partial_x v=&\frac{V(\rho)-v}{\tau}. \label{v}
\end{align}
The state variable $\rho(x,t)$ is the traffic density and $v(x,t)$ is the traffic speed. The equilibrium velocity-density relationship $V(\rho)$ is a decreasing function of density. The equilibrium flux function $Q(\rho)$ is 
\begin{align}
Q(\rho) = \rho V(\rho).
\end{align}
For ARZ model, the choice of $V(\rho)$ needs to satisfy that the flux function $Q(\rho)$ is strictly concave $Q(\rho)'' < 0$ for model validity and prediction of traffic waves. We choose the following model in the form of Greenshield's model,
\begin{equation}
V(\rho)=v_f\left(1-\left(\frac{\rho}{\rho_m}\right)^\gamma\right). \label{vf}
\end{equation}
Note that the choice of $V(\rho)$ here is due to its simplicity of expression and smoothness. The observer design proposed in this paper is not limited by this choice. 

The inhomogeneous ARZ including a relaxation term on the right hand side of the velocity PDE is considered. The constant parameter $\tau$ is the relaxation time which describes drivers' driving behavior adapting to equilibrium density-velocity relation over time. Note that the ARZ model without the relaxation term cannot address this phenomenon and poses an easier estimation problem. 

The increasing function of density $p(\rho)$ is defined as the traffic pressure
\begin{align}
p(\rho)=C_0\rho^\gamma,
\end{align} 
where $C_0,\gamma \in \mathbb{R}_+$.
The pressure function $p(\rho)$ is chosen so that it is related to equilibrium velocity-density function $V(\rho)$ as 
	\begin{align}
	p(\rho) = V(0) - V(\rho).
	\end{align}
Given $V(\rho)$ in \eqref{vf}, we have density pressure as
\begin{align}
	p(\rho) = v_f\left( \frac{\rho}{\rho_m^\gamma}\right)^\gamma.
\end{align}
Note that this choice of model parameter shows a marginal stability and a very slow damping effect of stop-and-go traffic. The following boundary observer design can be applied to the model when the above relation does not hold.

\subsection{Linearized ARZ model in traffic flux and velocity}
The traffic density is defined as the number of vehicles per unit length while the traffic flux represents the number of vehicles per unit time which cross a given point on the road. 
The traffic flow flux $q$ is 
\begin{align}
q=\rho v.
\end{align}
Traffic flux $q$ and velocity $v$ are most accessible physical variables to measure in freeway traffic. $q$ is commonly measured by loop detectors and $v$ is obtained by GPS. Therefore, we rewrite the ARZ model in traffic flux $q$ and traffic velocity $v$ as follows,
\begin{align}
q_t + vq_x =&\frac{q(\gamma  p-v)}{v} v_x+\frac{q(v_f-p-v)}{\tau v}, \label{nqv1}\\
v_t-(\gamma  p-v) v_x=&\frac{v_f-p-v}{\tau}, \label{nqv2}
\end{align}
$p(\rho)$ and $q$ are related by
\begin{align}
p=\frac{v_f}{\rho_m^\gamma}\left(\frac{q}{v}\right)^\gamma. \label{p}
\end{align} 
There is no explicit solution to the above nonlinear coupled hyperbolic system. To further understand the dynamics of ARZ traffic model in $(q,v)$-system, we linearize the model around steady states ($q^\star$, $v^\star$) which are chosen as spatial and temporal nominal values of state variables. Small deviations from the nominal profile are defined as
\begin{align}
\tilde q(x,t)=&q(x,t)-q^\star,\\
\tilde v(x,t)=&v(x,t)-v^\star.
\end{align}
The steady density is given as $\rho^\star =q^\star/v^\star $ and setpoint density-velocity relation satisfy the equilibrium relation $V(\rho)$,
\begin{align}
 v^\star = V(\rho^\star). 
\end{align}


The linearized ARZ model in $(\tilde q, \tilde v)$ around reference system $(q^\star, v^\star)$ with boundary conditions is given by
\begin{small}
\begin{align}
\notag\tilde{q}_t+v^\star\tilde{q}_x+\frac{q^\star}{v^\star}\left(v^\star + \frac{q^\star}{v^\star}V'\left(\frac{q^\star}{v^\star}\right)\right) \tilde{v}_x=&-{q^\star}\frac{(v^\star)^2+q^\star V'\left(\frac{q^\star}{v^\star}\right)}{\tau (v^\star)^3}\tilde v \\&+\frac{q^\star V'\left(\frac{q^\star}{v^\star}\right)}{\tau (v^\star)^2}\tilde q, \label{tq} \\
\notag \tilde v_t +\left( v^\star + \frac{q^\star}{v^\star}V'\left(\frac{q^\star}{v^\star}\right)\right) \tilde v_x =& -\frac{(v^\star)^2+q^\star V'\left(\frac{q^\star}{v^\star}\right)}{\tau (v^\star)^2}{\tilde v}\\&+\frac{V'\left(\frac{q^\star}{v^\star}\right)}{\tau v^\star}\tilde q, \label{tv}
\end{align}	
\end{small}
where the two characteristic speeds of the above linearized PDE model are
\begin{align}
\lambda_1 = &v^\star, \\ \lambda_2 =& v^\star + \rho^\star V'(\rho^\star).
\end{align}
\begin{itemize}
	\item Free-flow regime  $: \lambda_1> 0  , \; \lambda_2 >0$\\
	In the free-flow regime, both the disturbances of traffic flux and velocity travel downstream, at respective characteristic speeds $\lambda_1$ and $\lambda_2 $. The linearized ARZ model in free-regime is a homo-directional hyperbolic system.
	\\
	\item Congested regime $: \lambda_1> 0  , \; \lambda_2 <0$\\ In the congested regime, the traffic density is greater than a critical value and second characteristic speed becomes negative value. Therefore, disturbances of the traffic speed travel upstream with $\lambda_2$  while the disturbances of the traffic flow flux are carried downstream with the characteristic speed $\lambda_1$. The hetero-directional propagations of disturbances force vehicles into stop-and-go traffic. 
\end{itemize}
Since in the free-flow regime, the linearized homo-directional hyperbolic PDEs can be solved explicitly by the inlet boundary values and therefore state estimates can be obtained by solving the hyperbolic PDEs. In this work, we focus on the congested regime with two hetero-directional hyperbolic PDEs. It is a more relevant and challenging problem for traffic states estimation.

\section{Boundary Observer Design}
In this section, boundary sensing is employed for observer design. The state estimation of nonlinear ARZ model is achieved using backstepping method. The output injection gains are designed for the linearized ARZ model and then adding to the copy of nonlinear plant.

Boundary values of state variations from steady states are defined as 
\begin{align}
Y_{q,in}(t)  = & \tilde q(0,t), \label{y1}\\
Y_{q,out}(t)  = & \tilde q(L,t),\\
Y_v(t)  = & \tilde v(L,t).\label{y2}
\end{align}
where the values of $\tilde q(0,t)$, $\tilde q(L,t)$ and $\tilde v(L,t)$ are obtained by subtracting setpoint values from the sensing of incoming traffic flux $q(0,t)$, outgoing flux $q(L,t)$ and outgoing velocity $v(L,t)$,
\begin{align}
y_q(t) =& q(0,t), \\
y_{out}(t)= & q(L,t),\\
y_v(t) = & v(L,t).
\end{align}
Real-time sensing of traffic flux and velocity can be obtained by high-speed camera or induction loops. The induction loops are coils of wire embedded in the surface of the road to detect changes of inductance when vehicles pass. 

\subsection{Output injection for linearized ARZ model}
We diagonalize the linearized equations and therefore write $(\tilde q, \tilde v)$-system in the Riemann coordinates. The Riemann variables are defined as
\begin{align}
\xi_1 = &\frac{\rho^\star \lambda_2}{\lambda_1-\lambda_2}\tilde v + \tilde q, \\ \xi_2 =& \frac{q^\star}{\lambda_1-\lambda_2}\tilde v，
\end{align}
The inverse transformation is given by 
\begin{align}
\tilde v =& \frac{\lambda_{1} - \lambda_{2}}{q^\star} \xi_2,\\
\tilde q =& \xi_1 - \frac{\lambda_{2}}{\lambda_{1}}\xi_2.
\end{align}
The measurements are taken at boundaries result in the following boundary conditions
\begin{align}
\xi_1(0,t) =& \frac{\lambda_2}{\lambda_1} \xi_2(0,t) + Y_q(t),\\
\xi_2(L,t) = &  \frac{q^\star}{\lambda_1-\lambda_2} Y_v(t).
\end{align}
Therefore the linearized ARZ model in Riemann coordinates is obtained
\begin{align}
\partial_t \xi_1 + \lambda_1 \partial_x \xi_1 =& -\frac{1}{\tau}\xi_1,\\
\partial_t \xi_2 + \lambda_2 \partial_x \xi_2 =& -\frac{1}{\tau}\xi_1,\\
\xi_1(0,t) =& \frac{\lambda_2}{\lambda_1} \xi_2(0,t),\\
\xi_2(L,t) = & \xi_1(L,t).
\end{align}

In order to diagonalize the right hand side to implement backstepping method, we introduce a scaled state as follows:
\begin{align}
\bar w(x,t) = & \exp\left(\frac{ x }{\tau \lambda_1 } \right)  \xi_1(x,t), \label{tr1}\\
\bar v(x,t) = & \xi_2(x,t).
\end{align}

The $(\xi_1, \xi_2 )$-system is then transformed to a first-order $2 \times 2$ hyperbolic system
\begin{align}
\bar w_t(x,t)  + \lambda_1 \bar w_x(x,t) =& 0, \label{bw}\\
\bar v_t(x,t) + \lambda_2 \bar v_x(x,t) =& c(x) \bar w(x,t), \label{bv}\\
\bar w(0,t)=&\frac{\lambda_2}{\lambda_1} \bar v(0,t) + Y_{q,in}(t) ,\label{bbvw0}\\
\bar v(L,t)=& \frac{q^\star}{\lambda_1-\lambda_2} Y_v(t), \label{bbvw1}
\end{align}
where the spatially varying parameter $c(x)$ is defined as 
\begin{align}
c(x)=-\frac{1}{\tau}\exp\left(-\frac{ x }{\tau \lambda_{1} } \right),
\end{align}	
Parameter $c(x)$ is a strictly increasing function and bounded by
\begin{align}
-\frac{1}{\tau}\leq c(x)\leq -\frac{1}{\tau}\exp\left(-\frac{L}{\tau \lambda_{1} } \right) . \label{cbound}
\end{align}

Then we design a boundary observer for the linearized ARZ model to estimate $\bar w(x,t)$ and $\bar v(x,t)$ by constructing the following system
\begin{align}
\hat w_t(x,t)  + \lambda_1 \hat w_x(x,t) =& r(x)(\bar w(L,t)-\hat w(L,t))  ,\label{O1}\\
\notag \hat v_t(x,t) + \lambda_2 \hat v_x(x,t) =& c(x) \hat w(x,t)\\
&+ s(x)(\bar w(L,t)-\hat w(L,t))  , \\ 
\hat w(0,t)=&
\frac{\lambda_2}{\lambda_1} \hat v(0,t) + Y_{q,in}(t),\\
\hat v(L,t)
=&\frac{q^\star}{\lambda_1-\lambda_2} Y_v(t) , \label{2}
\end{align}
where $\hat w(x,t)$ and $\hat v(x,t)$ are the estimates of the state variables $\bar w(x,t)$ and $\bar v(x,t)$. The value $\bar w(L,t)$ is obtained by plugging in the measured outgoing flow flux $Y_{q,out}(t)$ and velocity $Y_v(t)$ into \eqref{tr1},
\begin{align}
\bar w(L,t) = \exp\left(\frac{ L }{\tau \lambda_1 } \right)\left(\frac{\rho^\star \lambda_2}{\lambda_1-\lambda_2} Y_v(t) + Y_{q,out}(t)\right). \label{bwl}
\end{align}

The term $r(x)$ and $s(x)$ are output injection gains to be designed. We denote estimation errors as
\begin{align}
	\check w(x,t)=&\bar w(x,t) -\hat w(x,t),\\
	\check v(x,t)=&\bar v(x,t) -\hat v(x,t).
\end{align}
The error system is obtained by subtracting the estimates \eqref{O1}-\eqref{2} from \eqref{bw}-\eqref{bbvw1},
\begin{align}
\check w_t(x,t) =& -v^\star \check w_x(x,t)-r(x)\check w(L,t), \label{e21}\\
\notag	\check v_t(x,t) =& (\gamma p^\star - v^\star )\check v_x(x,t)+c(x) \check w(x,t)\\
&-s(x)\check w(L,t),\\ 
\check w(0,t)=&
\frac{\lambda_2}{\lambda_1} \check v(0,t), \\
\check v(L,t)
=& 0.\label{e24}
\end{align}

The output injection gains $r(x)$ and $s(x)$ need to guarantee the error system $(\check w, \check v)$ decays to zero. Using backstepping transformation, we transform the error system \eqref{e21}-\eqref{e24} into the following target system
\begin{align}
\alpha_t(x,t) + \lambda_1 \alpha_x(x,t) =& 0, \label{al}\\
\beta_t(x,t) + \lambda_2 \beta_x(x,t) =& 0,\\ 
\alpha (0,t)=&\frac{\lambda_2}{\lambda_1} \beta(0,t),  \label{inb}\\
\beta(L,t)=& 0.\label{outb}
\end{align}
The explicit solution to the target system \eqref{al}-\eqref{outb} can be easily found. Thus we have 
\begin{align}
\alpha(x,t)\equiv \beta(x,t)\equiv 0,
\end{align}
after a finite time $t=t_f$ where
\begin{align}
t_f = \frac{L}{|\lambda_{1}|} + \frac{L}{|\lambda_{2}|}. \label{tf}
\end{align}
It is straightforward to prove that the $\alpha, \beta$ system is $L^2$ exponentially stable. 

The backstepping transformation is
\begin{align}
\alpha(x,t) =& \check w(x,t)-\int^{L}_{x} K(L+x-\xi) \check w(\xi,t)d\xi,\\
\beta(x,t)=& \check v(x,t)-\int^{L}_{x} M(\lambda_1  x - \lambda_2\xi)\check w(\xi,t)d\xi,
\end{align}
where the kernel variables $K(x)$ and $M(x)$ map the error system into the target system. The kernel $M(x)$ is defined as
\begin{align}
M(x)=-\frac{1}{\lambda_1  - \lambda_2} c\left(\frac{x}{\lambda_1 - \lambda_2}\right).
\end{align}
For boundary condition \eqref{inb} to hold, the kernels $ K(x)$ and $\check M(x)$ satisfy the relation
\begin{align}
K(L-\xi)=& M ((\lambda_2  - \lambda_1)\xi).
\end{align}
the kernel $K$ is then obtained
\begin{align}
K(x)=-\frac{1}{\lambda_1  - \lambda_2} c\left(\frac{-\lambda_2}{\lambda_1 - \lambda_2}(L-x)\right).
\end{align}
According to the boundedness of $c(x)$ in \eqref{cbound}, the kernels are bounded
\begin{align}
|K(x)| \leq \frac{1}{(\lambda_1 - \lambda_2) \tau}, \label{bck}
\end{align}
and therefore $M(x)$ is bounded.
The output injection gain $r(x)$ and $s(x)$ are given by
\begin{align}
r(x)=&\lambda_1  K(x)=-\frac{\lambda_1}{\lambda_1 - \lambda_2} c\left(\frac{ \lambda_2}{\lambda_1 - \lambda_2}(L-x)\right), \label{g1}\\
\notag	s(x)=&- \lambda_1 M( \lambda_1  x - \lambda_2 L)\\
=&\frac{\lambda_1}{\lambda_1 - \lambda_2}  c\left( x - \frac{ \lambda_2}{\lambda_1 - \lambda_2}(L-x)\right).\label{g2}
\end{align}
The backstepping transformation is invertible. Therefore, we study the stability of the error system through the target system \eqref{al}-\eqref{outb}. It is straightforward to prove the exponential stability of error system in the $L^2$ sense and finite-time convergence in $t_f$ given by \eqref{tf}. We arrive at the following theorem.

\begin{thm}
	Consider system \eqref{e21}-\eqref{e24} with inital conditions $\check w_0, \check v_0 \in L^2[0,L]$. The equilibrium $\check w \equiv \check v \equiv 0$ is exponentially stable in the $L^2$ sense. It holds that
	\begin{align}
		||\bar w(\cdot,t)-\hat w(\cdot, t)||\to 0\\
		||\bar v(\cdot,t)-\hat v(\cdot, t)||\to 0
	\end{align}
	 and the convergence to equilibrium is reached in finite time $t=t_f$ given in \eqref{tf}.
\end{thm}

\subsection{Boundary observer design for Nonlinear ARZ model}
For nonlinear boundary observer, we construct the system by keeping the output injection that is designed for the linearized ARZ model, then add them to the copy of original nonlinear ARZ model.

We summarize the transformation from linearized ARZ model in $(\tilde q,\tilde v)$-system to $(\bar w, \bar v)$-system,
\begin{align}
\bar w(x,t)=& \exp\left(\frac{ x }{\tau \lambda_1 } \right)\left(\frac{\rho^\star \lambda_2}{\lambda_1-\lambda_2}\tilde v + \tilde q\right),\label{tf1}\\
\bar v(x,t)=&\frac{q^\star}{\lambda_1-\lambda_2}\tilde{v}(\xi,t).\label{tf2}
\end{align}
And the inverse transformation is given by
\begin{align}
\tilde q(x,t)=&\exp\left(-\frac{ x }{\tau  \lambda_1 } \right)\bar w(x,t) -  \frac{\lambda_2}{ \lambda_1}\bar v(x,t), \label{qv1}\\
\tilde v(x,t)=&\frac{\lambda_1-\lambda_2}{q^\star} \bar{v}(\xi,t) .\label{qv2}
\end{align}
Due to equivalence of $(\check w,\check v)$ and $(\tilde q, \tilde v)$-system, we arrive at the following theorem for the linearized ARZ model.
\begin{thm}
	Consider system \eqref{tq}-\eqref{tv} with inital conditions $\tilde q_0, \tilde v_0 \in L^2[0,L]$. The equilibrium $\tilde q \equiv \tilde v \equiv 0$ is exponentially stable in the $L^2$ sense. It holds that
	\begin{align}
	||q (\cdot,t)-q^\star||\to 0\\
	||v (\cdot,t)-v^\star||\to 0
	\end{align}
	and the convergence to set points is reached in finite time $t=t_f$.
\end{thm}

We denote the error injections designed for the linearized ARZ model \eqref{O1}-\eqref{2} as
\begin{align}
E_w(t) =& -r(x)(\bar w(L,t)-\hat w(L,t)), \\
E_v(t) =& -s(x)(\bar w(L,t)-\hat w(L,t)).
\end{align}
The output injection gains $r(x)$, $s(x)$ are designed in \eqref{g1} and \eqref{g2}. According to \eqref{bwl}, $\bar w(L,t)$ is obtained from real-time measurement of the traffic boundary data in \eqref{y1}-\eqref{y2}. Therefore, the values of output injections $E_w(t)$ and $E_v(t)$ are known. 

The nonlinear observer $(\hat\rho(x,t),\hat v(x,t))$ for state estimation is obtained by combining the copy of the nonlinear ARZ model $(\rho,v)$ given by \eqref{rho}\eqref{v} and the above linear injection errors in original state variables density and velocity,
\begin{align}
\partial_t \hat \rho + \partial_x( \hat \rho \hat v)=&  { \frac{1}{v^\star}\left(\exp\left(-\frac{ L }{\tau  \lambda_1 } \right)E_w- E_v\right)}, \\
\partial_t \hat v+(\hat v + \hat \rho V'(\hat \rho))\partial_x \hat v=&\frac{V(\hat \rho)- \hat v}{\tau} +{ \frac{\lambda_1-\lambda_2}{q^\star}  E_v},
\end{align}
where the linear injection on the right hand side are obtained by transforming $(\hat w, \hat v)$ to $(\rho, v)$ given in \eqref{qv1},\eqref{qv2}. The boundary conditions are
\begin{align}
\hat \rho(0,t) = \frac{y_q(t)}{\hat v(0,t)},\\
\hat v(L,t) = {y_v(t)}.
\end{align}
When the initial states of the system is close to the set points, the linearized part dominates the nonlinear estimation error system and therefore we can prove that the local $H^2$ exponential stability holds for estimation error system of the nonlinear ARZ model,following approach in \cite{rafael}. This observer result is validated in the following numerical simulation.

\begin{table}[t!]
	\caption{Parameter Table}
	\label{tab:default}
	\centering
	\begin{tabular}{ p{5cm}|p{3cm}  }	
		\hline
		Parameter Name  & Value \\
		\hline
		Maximum traffic density $\rho_m$   & 160 vehicles/km \\
		Traffic pressure and coefficient $\gamma$ & 1 \\
		Maximum traffic velocity $v_f$	&  40 m/s  \\
		Relaxation time $\tau$  & 60 s \\
		Reference density $\rho^\star$ & 120 vehicles/km \\
		Reference velocity $v^\star$ & 10 m/s\\
		Freeway segment length $L$ & $500$ m
		\\		\hline
	\end{tabular}
\end{table}

\section{Simulation}

For simulation of nonlinear ARZ PDE model, we assume that
the initial conditions are sinusoidal oscillations around the steady states. Model parameters are chosen as shown in the table 1. We consider a constant incoming flow and constant outgoing density for boundary conditions,
\begin{align}
	\tilde q(0,t) = & 0, \label{bc1}\\
	\tilde v(L,t) = & \frac{1}{\rho^\star}\tilde q(L,t).\label{bc2}	
\end{align}
When dealing with the real traffic data, we do not prescribe any boundary conditions beforehand but import the boundary data directly.

We use finite volume method and write the ARZ model in the conservative variables, then apply two-stage Lax-Wendroff scheme to discretize ARZ model in spatio-temporal domain. The scheme is second-order accurate in space. The spatial grid resolution is chosen to be smaller than the average vehicle size so that the numerical errors are smaller than the model errors. Therefore the simulation is valid for this continuum model. 


For the numerical stability of Lax - Wendroff scheme, the spatial grid size $\Delta x$ and time step $\Delta t$ is chosen so that CFL condition is satisfied:
\begin{align}
\max |\lambda_{1,2}| \leq \frac{\Delta x}{\Delta t},
\end{align} 
%

Note that we need to specify state values at both $x=0$ and $x=L$ boundaries. ARZ model will pick up some combination of $\rho$ and $v$ at each of the two boundaries, depending on the direction of characteristics. We implement the boundary conditions in \eqref{bc1} and \eqref{bc2}.


\begin{figure}[t!]
	\includegraphics[width=9.5cm]{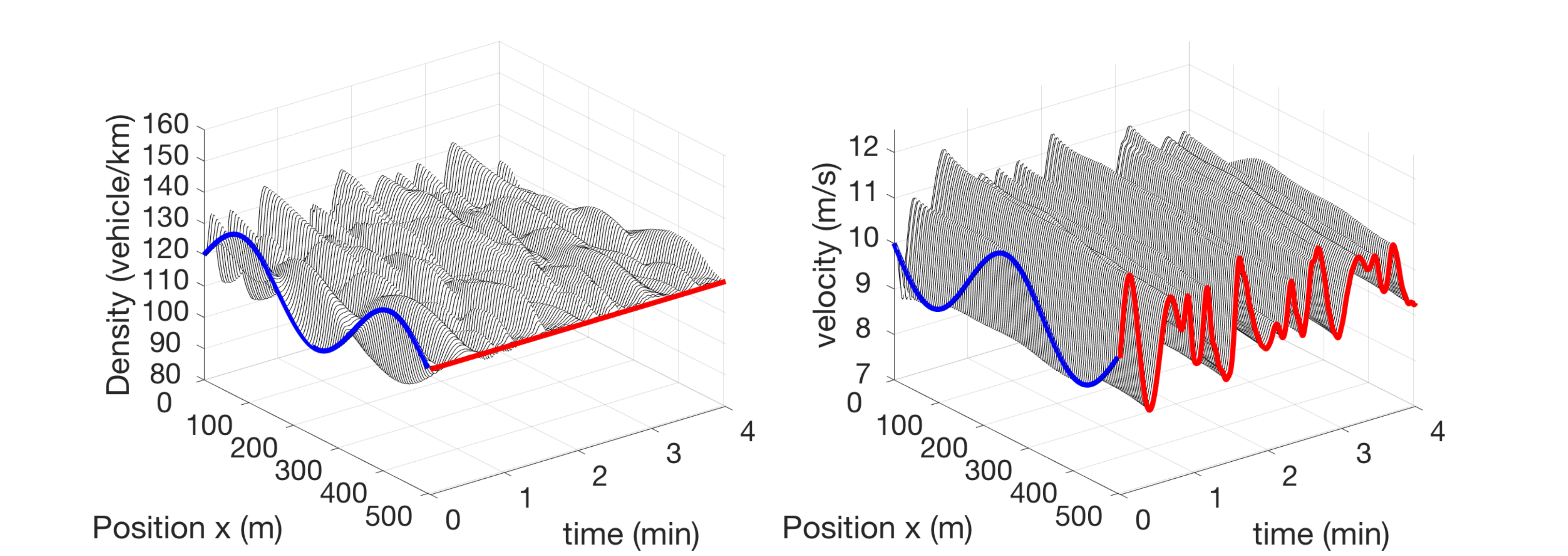}
	\centering
	\caption{Density $\rho(x,t)$ and velocity $v(x,t)$ of nonlinear ARZ model.}
\end{figure} 
\begin{figure}[t!]
	\includegraphics[width=9.5cm]{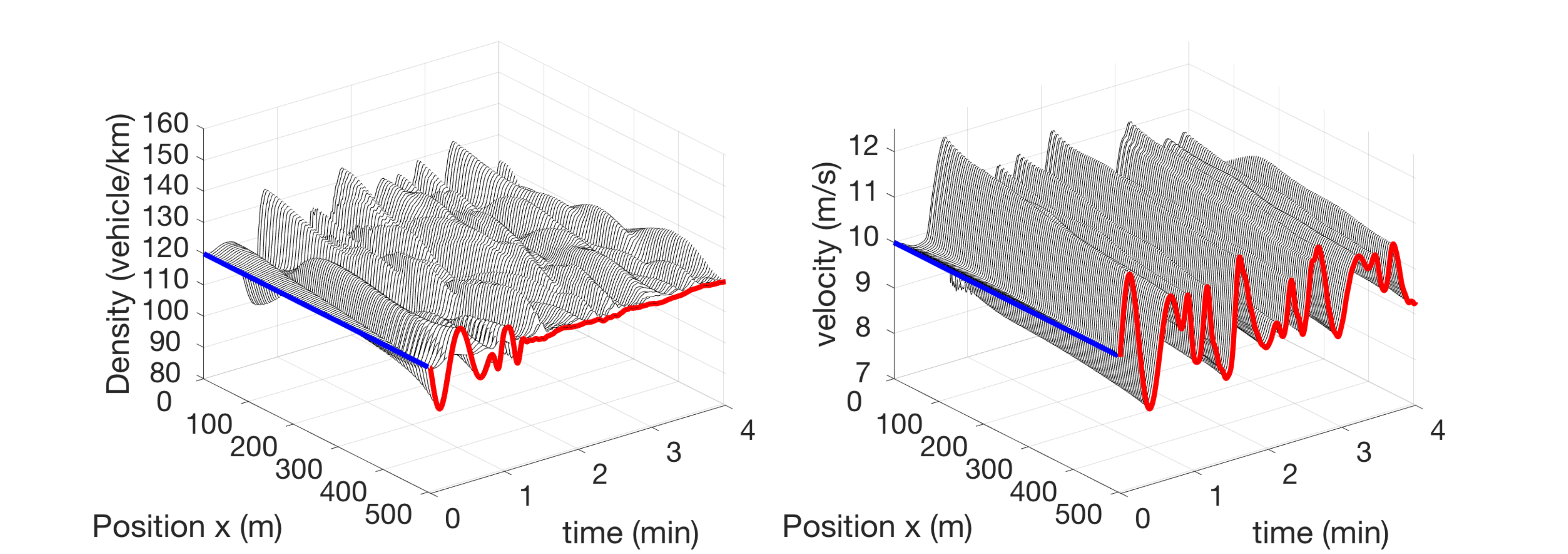}
	\centering
  \caption{States estimates $\hat \rho(x,t)$ and $\hat v(x,t)$ of nonlinear boundary observer.}
\end{figure} 
\begin{figure}[t!]
	\includegraphics[width=9.5cm]{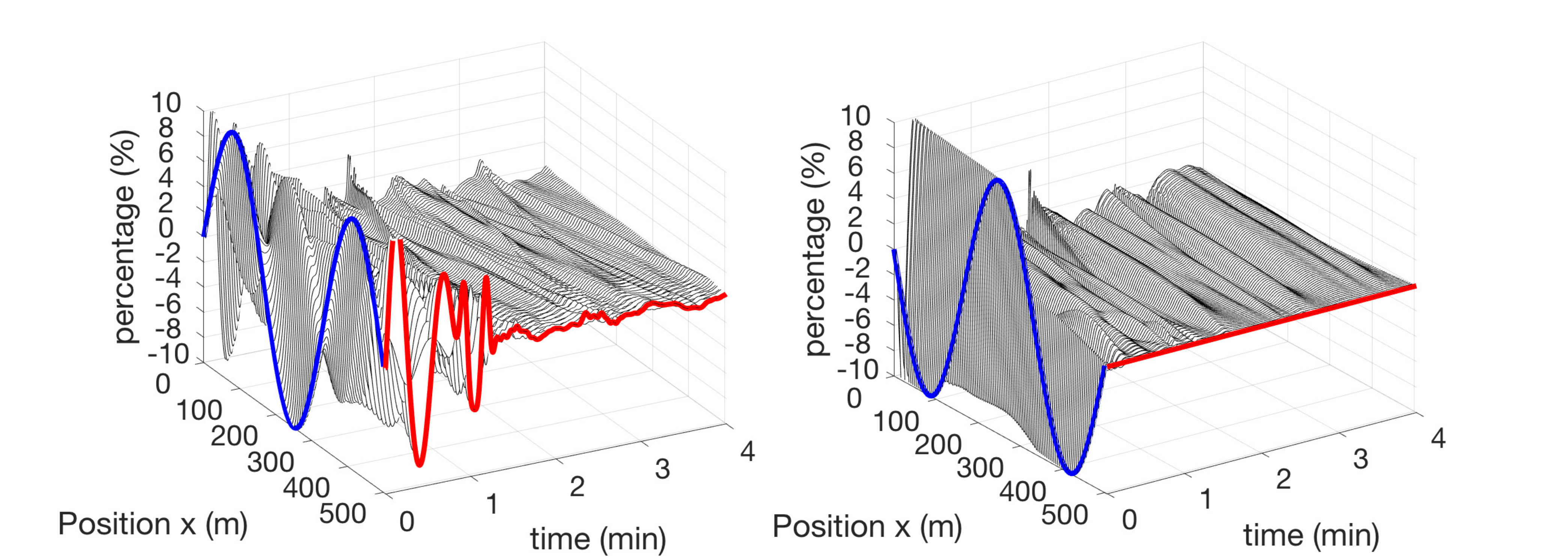}
	\centering
	\caption{Estimation errors $\check \rho(x,t)$ and $\check v(x,t)$.}
\end{figure}

In Fig 1-3, blue lines represent initial conditions while the red lines represent the evolution of outlet state value in temporal domain. The simulation is performed for a $500 \; \rm m$ length of freeway segment and evolution of traffic states density and velocity are plotted for $4 \; \rm min$. 

In Fig 1, traffic density and velocity are slightly damped and keeps oscillating. It takes the initial disturbance-generated vehicles to leave the domain in $50\;\rm s$ but the oscillations sustain for more than $4 \; \rm min$. Since the traffic states are chosen to be in the congested regime, the stop and go traffic appears in the simulation. 

State estimation of traffic density and velocity by nonlinear observer is shown in Fig.2. The measurement is taken for outgoing velocity and outgoing flow. The incoming flow is assumed to be at setpoint traffic flux. We do not assume any prior knowledge of initial condition and thus set the initial conditions to be at setpoint density and velocity. We can see that state estimates converges to true values of plant after $75 \; \rm s$. 

In Fig 3, the evolution of estimation errors are shown. After $75 \; \rm s$, the estimation errors for density and velocity converge to value less than $2\%$ of the setpoint value. There are still relatively very small estimation errors remain in the domain for two reasons. Our result only guarantees convergence of estimates in the spatial $L^2$ norm. In addition, there could be nonlinearities remain in the linear injection of nonlinear boundary observer design.


\section{Conclusion}

In conclusion, we develop a nonlinear boundary observer to estimate traffic states for the ARZ PDE model. Analysis of the linearized model leads our main focus to the congested regime where stop-and-go happens. Using spatial transformation and PDE backstepping method, we construct a boundary observer with a copy of the nonlinear plant and output injection of measurement errors so that local exponential stability of estimation errors in the $L^2$ norm and finite-time convergence to zeros are guaranteed. The observer design for nonlinear ARZ model has been validated by simulation for a stretch of freeway. For future work, validation of the nonlinear observer with traffic field data is of our interest.

\end{document}